\documentclass[conference]{IEEEtran}
\IEEEoverridecommandlockouts
\usepackage{cite}
\usepackage{amsmath,amssymb,amsfonts}
\usepackage{nccmath}
\usepackage{algorithmic}
\usepackage[ruled,vlined]{algorithm2e}
\usepackage{graphicx}
\usepackage{textcomp}
\usepackage{nomencl}
\makenomenclature
\usepackage{xcolor}
\def\BibTeX{{\rm B\kern-.05em{\sc i\kern-.025em b}\kern-.08em
    T\kern-.1667em\lower.7ex\hbox{E}\kern-.125emX}}
\begin{document}

\title{Designing Facilities to Improve Flexibility: Zone-based Dynamic Facility Layout with Embedded Input/Output Points
}

\author{\IEEEauthorblockN{Sadan Kulturel-Konak}
\IEEEauthorblockA{\textit{Management Information Systems} \\
\textit{Penn State Berks}\\
Reading, USA \\
sxk70@psu.edu}
\and
\IEEEauthorblockN{Abdullah Konak}
\IEEEauthorblockA{\textit{Information Sciences and Technology}\\
\textit{
Penn State Berks}\\
Reading, USA \\
auk3@psu.edu}
}

\maketitle

\begin{abstract}
This paper considers solving the unequal area Dynamic Facility Layout Problem (DFLP) using a zone-based structure. Zone-based layouts have significant advantages, such as being easily transferable to a detailed layout with innately included possible aisle structures; therefore, they can be fitted to the unique needs of the layout designers. The unequal area DFLP is modeled and solved using a zone-based structure, which is referred to as ZDFLP, where the dimensions of the departments and material handling system input/output (I/O) points are decision variables. A two-phase matheuristic, which directly operates on Problem ZDFLP without requiring an encoding scheme of the problem, is proposed to solve the ZDFLP with promising results. 
\end{abstract}


\section{Introduction}
Designing a facility over a multi-period planning horizon where the interdepartmental material flows change over the planning periods due to frequent changes in product demands forms a Dynamic Facility Layout Problem (DFLP). The DFLP, which was first introduced by Rosenblatt \cite{Rosenblatt1986}, in the continuous plane is a very challenging nonlinear optimization problem. This study considers a zone-based block layout to design manufacturing and logistics facilities considering material handling infrastructure. It is important to start with a low-cost block layout design which is the precursor to a detailed layout. Hence, it is essential to have a block layout that minimizes the material handling cost and is easily transferable to the related actual facility layout. After a flexible block layout has been accomplished, a detailed layout, including aisle structures, input/output points of the departments, and the exact locations of the departments, will be explored. A zone-based block layout \cite{Montreuil2002, Montreuil2004} inherently includes possible aisle structures which can quickly be adapted to different material handling systems and, therefore, can be transferred to a detailed layout with fewer modifications than a block layout based on the unrestricted general formulation that the Facility Layout Problem (FLP) would require. This is particularly important in the DFLP because the changes in a block layout from one period to the next may require structural modifications in the material handling system, which in turn may be very costly or, in some cases, impossible to implement practically. Currently, only a limited number of current models consider re-layout costs as a result of changes in the overall structure of the facility \cite{Xiao2017, Kulturel-Konak2019}.

Another point that needs to be addressed in the DLFP is the cost of re-purposing facility space for different department types. The relayout cost is generally evaluated as a function of the distance that departments are relocated during the redesign process. However, re-purposing a space unit from one department to another depends on the type of department. In healthcare facilities, for example, it may cost much less to re-purpose a regular treatment room as an intensive care unit than to re-purpose a waiting room. The cost of re-purposing space is particularly important in cyclic FLP \cite{Kulturel-Konak2015}, where the facility layout is periodically changed. This change is because of demand cycles or cases where the facility is temporarily repurposed to respond to the needs of catastrophic events such as pandemics or natural disasters and returned to its original functionality afterward. As the static version of the facility layout problem, Cubukcuoglu \cite{cubukcuoglu2022} suggested a hierarchical framework that divides the main design stages in hospitals into stacking (dividing the floors between functional spaces), zoning (placing rooms), and routing (embedding corridors).    

In this paper, we first present new formulations of the zone-based DFLP by considering the relayout cost of facility structural changes and space repurposing. Then, we present a matheuristic to solve this computationally difficult problem efficiently. Pérez-Gosende et al. \cite{Perez-Gosende2021}, in their DFLP review paper, also emphasized an emerging need for developing and applying more powerful matheuristic approaches as solution strategies to those models and integrating the economic, environmental, and social sustainable aspects into DFLP models.
 
\section{Background}

In general, layouts can be planned for brand-new plants, i.e., greenfield layout design, or existing plants, i.e., re-layout. The DFLP can be considered to combine both greenfield and re-layout design aspects. In the literature, more attention has been paid to greenfield designs, where existing restrictions have not influenced the layout plan. Although its limited importance/appearance in the literature, the re-layout problem is more frequent in practice \cite{Kulturel-Konak2007}. Recently, Pérez-Gosende et al. \cite{Perez-Gosende2021} reviewed articles from the layout literature based on the problem type, planning stages, material handling configurations, and solution approaches. Of these articles reviewed by Pérez-Gosende et al. \cite{Perez-Gosende2021}, only 11.21 percent dealt with the re-layout aspect of the layout problems. In a DFLP review article \cite{Zhu2018}, metaheuristics and hybrid approaches were mentioned as the most promising approaches for tackling complicated and realistic scenarios.
This paper uses a zone-based layout \cite{Montreuil2002, Montreuil2004} to solve the DFLP and offer a design approach to increase flexibility in layout design \cite{Xiao2017, Kulturel-Konak2019}. Zone-based layouts have significant advantages, such as they are easy to transfer to a detailed layout with innately included possible aisle structures, and, therefore, they can be fitted to the unique needs of the layout designers. Thus, zone-based designs are pertinent to many real-life production systems \cite{Montreuil2002}. These advantages are particularly significant in the DFLP because the changes in a block layout between consecutive periods may require structural modifications in the material handling system or the facility, which may be very costly or sometimes impractical to implement. For example, a zone-based layout that inherently includes possible aisle structures might be adapted in a recent bottom-up multi-objective DFLP model suggested by Pérez-Gosende et al. \cite{perez2023} while integrating corridors into the layout suggested by the authors as a future research area.

\section{Modelling Approach}
\begin{figure}
\centerline{\includegraphics{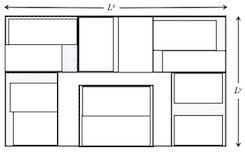}}
\caption{A sample zone-based layout with six zones. (Within a zone, the departments are allowed to be located vertically or horizontally.)}
\label{fig}
\end{figure}

The unequal area DFLP is modeled and solved using a zone-based structure where the dimensions of the departments and material handling system input/output (I/O) points are decision variables, unlike many previous DFLP models. A zone is defined as a sub-region in the facility with flexible locations and boundaries. Fig.~\ref{fig} shows a sample zone-based layout with six zones arranged in a U-shape layout. Although the zone locations and shapes are decision variables in the zone-based DFLP (ZDFLP) defined in this paper, facility designers may prefer defining the relative locations of the zones in the facility, which can be implemented in the ZDFLP model by fixing the decision variables to set zone relative locations.   

A zone-based block layout inherently defines aisle structures, which can be adapted to the current material handling system. In the model, the departments are arranged in flexible zones in either the $x$-axis or $y$-axis directions, but not both. In \cite{Xiao2017, Kulturel-Konak2019}, the departments are allowed to be placed freely within zones. However, this flexibility may lead to block layouts in which departments are located far from zone boundaries. On the other hand, arranging departments in either the $x$-axis or $y$-axis directions within a zone ensures that each department has access to a zone boundary for locating its I/O point. In other words, the ZDFLP model defined in this paper assumes that the department I/O points will be located adjacent to the aisles, which is well justified in real-life. In fact, the flexible bay structure, which is a form of zone-based structure, is preferred as a block layout schema in many cases because it can allow the designing of aisles that can be easily reconfigured or adapted to accommodate different production processes (see \cite{norman2001,kulturel2011} for more details about the advantages of the flexible bay structure.)  

In the formulation, the boundaries of the zones are flexible and can change from one period to another. Such structural changes are penalized in the objective function since they may require updating the material handling system and department interface points.

The notation used in the definition of the ZDFLP is given below. Fig. \ref{model} provides the MIP model without the linearization of the constraints with the absolute value function for the brevity of the model presentation.

\renewcommand{\nomgroup}[1]{%
\ifthenelse{\equal{#1}{V}}{\item[\textbf{Decision Variables\\}]}{%
\ifthenelse{\equal{#1}{P}}{\item[\textbf{Parameters}]}{}}}
\nomenclature[P,00]{$R$}{set of axis directions, $R=\{x,y\}$}
\nomenclature[P,00]{$S$}{set of cardinal directions, (e)ast, (w)est, (s)outh, and (n)orth, $S=\{e,w,s,n\}$}
\nomenclature[P,01]{$T$}{set of planning periods.}
\nomenclature[P,03]{$\Omega_{t}$}{set of departments used in period $t$.}
\nomenclature[P,05]{$K$}{set of zones in the facility.}
\nomenclature[P,11]{$f_{ijt}$}{material flow between departments $i$ and $j$ in period $t$.}
\nomenclature[P,13]{$L^r$}{side length of the facility in the $r$-axis direction.}
\nomenclature[P,15]{$a_{it}$}{minimum area requirement of department $i$ in period $t$.}
\nomenclature[P,18]{$\underline{l}_{it}^r$}{minimum allowed side length of department $i$ in the $r$-axis direction in period $t$.}
\nomenclature[P,21]{$\bar{l}_{it}^r$}{maximum allowed side length of department $i$ in the $r$-axis direction in period $t$.}
\nomenclature[P,24]{$C_{ijt}$}{cost of moving one unit material per unit distance between departments $i$ and $j$ in period $t$.}
\nomenclature[P,25]{$TC$}{the total cost of the layouts in the planning horizon.}
\nomenclature[P,27]{$R_{it}$}{fixed rearrangement cost of department $i$ between periods $t-1$ and $t$.}
\nomenclature[P,30]{$Q_{it}$}{unit distance rearrangement cost of department $i$ between periods $t-1$ and $t$.}
\nomenclature[P,32]{$B_{kt}$}{fixed cost of moving the boundaries of zone $k$ between periods $t-1$ and $t$.}
\nomenclature[P,35]{$\Delta$}{number of tangential support points to approximate department areas.}
\nomenclature[P,36]{$\bar{x}_{ipt}$}{tangential support point $p$ for department $i$ in period $t$.}

\nomenclature[V,00]{$\beta_k$}{binary variable indicating zone type such that departments are aligned along the $x$-axis direction ($\beta_k=0$) or $y$-axis direction ($\beta_k=1$).}
\nomenclature[V,01]{$z_{ijt}^r$}{the relative locations of departments $i$ and $j$ in period $t$ such that $z_{ijt}^r=1$  if department $i$ is enforced to precede department $j$ in the $r$-axis direction, $z_{ijt}^r=0$, otherwise.}
\nomenclature[V,02]{$\gamma_{kh}^r$}{the relative locations of zones $k$ and $h$ in period $t$ such that $\gamma_{kh}^r=1$  if zone $k$ is enforced to precede zone $h$ in the $r$-axis direction, $\gamma_{kh}^r=0$, otherwise.}
\nomenclature[V,021]{$b_{ikt}$}{department to zone assignment such that $b_{ikt}=1$ if department $i$ is assigned to zone $k$ in period $t$, $b_{ikt}=0$, otherwise.}
\nomenclature[V,03]{$l_{it}^{r}$}{half side length of department $i$ in the $r$-axis direction in period $t$.}
\nomenclature[V,05]{$c_{it}^{r}$}{center coordinates of department $i$ in the $r$-axis direction in period $t$.}	
\nomenclature[V,06]{$u_{it}^{r}$}{amount of the change of the center of department $i$ in the $r$-axis direction between periods $t-1$ and $t$.}
\nomenclature[V,07]{$g_{it}^{r}$}{I/O point coordinates of department $i$ in the $r$-axis direction in period $t$.}	 
\nomenclature[V,08]{$d_{ijt}^{r}$}{distances between the input/output points of departments $i$ and $j$ in the $r$-axis direction in period $t$.}	 
\nomenclature[V,11]{$v_{it}$}{department rearrangements such that $v_{it}$=1 if department $i$ is relocated between periods $t-1$ and $t$, $v_{it}=0$ otherwise.}
\nomenclature[V,13]{$q_{it}^{r}$}{the coordinate of $r \in S$ bound of zone $k$ in period $t$.}
\nomenclature[V,15]{$o_{kt}^r$}{zone boundary arrangement such that $o_{kt}^r=1$ if $q_{it}^{r}$ is relocated between periods $t-1$ and $t$, $o_{kt}^r=0$ otherwise.}
\renewcommand{\nomname}{}
\printnomenclature

\begin{figure*}[!htb]
\begin{align}
min \quad TC=&\sum_{t \in T} \sum_{(i,j) \in P_t} \sum_{r \in R} C_{ijt}f_{ijt}d_{ijt}^r + && \label{eq:modelo1} \\ 
&\sum_{t \in T: t>1} \sum_{i \in \Omega_{t}} R_{it} v_{it} + \sum_{t \in T: t>1} \sum_{i \in \Omega_{t}} \sum_{r \in R}Q_{it}u_{it}^r+&& \label{eq:modelo2} \\
&\sum_{t \in T:t>1}\sum_{k \in K}\sum_{r \in S} B_{kt}o_{kt}^r& & \label{eq:modelo3} \\
s.t. & & & \nonumber\\
&\gamma_{kht}^x+\gamma_{hkt}^x+\gamma_{kht}^y+\gamma_{hkt}^y = 1 & \forall t \in T, \{k,h\} \in K : k < h & \label{eq:modelczone} \\
&q_{kt}^w \leq q_{kt}^e  &  \forall t \in T, k \in K & \label{eq:modelcz7} \\
&q_{kt}^e \leq q_{ht}^w + L^x(1-\gamma_{kht}^x)  &  \forall t \in T, \{k,h\} \in K: k \neq h & \label{eq:modelcz8} \\
&q_{kt}^s \leq q_{kt}^n  &  \forall t \in T, k \in K  & \label{eq:modelcz9} \\
&q_{kt}^n \leq q_{ht}^s + L^y(1-\gamma_{kht}^y) &  \forall t \in T, \{k,h\} \in K: k \neq h & \label{eq:modelcz10} \\
&q_{kt}^e \leq L^x  &  \forall t \in T, k \in K  & \label{eq:modelcz11} \\
&q_{kt}^n \leq L^y  &  \forall t \in T, k \in K  & \label{eq:modelcz12} \\
&q_{kt}^w \geq 0  & \forall t \in T, k \in K  & \label{eq:modelcz13} \\
&q_{kt}^s \geq 0  & \forall t \in T, k \in K  & \label{eq:modelcz14} \\
&z_{ijt}^x+z_{jit}^x \geq b_{ikt}+b_{jkt} -1-\beta_k & \forall t \in T, k \in K, \{i, j\} \in \Omega_{t} : i < j & \label{eq:modelc1} \\
&z_{ijt}^y+z_{jit}^y \geq b_{ikt}+b_{jkt} -2+\beta_k & \forall t \in T, k \in K, \{i, j\} \in \Omega_{t} : i < j & \label{eq:modelc2} \\
&z_{ijt}^r+z_{jit}^r \leq 1 &  \forall t \in T, r \in R, \{i,j\} \in \Omega_{t} : i < j  & \label{eq:modelc3} \\
&c_{it}^r+l_{it}^r \leq c_{jt}^r-l_{jt}^r +L^r(1-z_{ijt}^r) & \forall t \in T, \{i,j\} \in \Omega_{t}, r \in R : i \neq j & \label{eq:modelc4} \\
&\sum_{k \in K}b_{ikt}=1 &  t \in T, i \in \Omega_{t} & \label{eq:model7} \\
&\sum_{i \in \Omega_{t}}b_{ikt} \geq 1 &  t \in T, k \in K  & \label{eq:model8} \\
&c_{it}^x+l_{it}^x \leq q_{kt}^e+ L^x(1-b_{ikt}) & \forall t \in T, k \in K, i \in \Omega_{t} & \label{eq:modelcz1} \\
&c_{it}^x-l_{it}^x \geq q_{kt}^w- L^x(1-b_{ikt}) & \forall t \in T, k \in K, i \in \Omega_{t} & \label{eq:modelcz2}\\
&c_{it}^y+l_{it}^y \leq q_{kt}^n+ L^y(1-b_{ikt}) &  \forall t \in T, k \in K, i \in \Omega_{t} & \label{eq:modelcz3}\\
&c_{it}^y-l_{it}^y \geq q_{kt}^s- L^y(1-b_{ikt}) &  \forall t \in T, k \in K, i \in \Omega_{t} & \label{eq:modelcz4}\\
&Mo_{kt}^r \geq |q_{kt}^r-q_{k(t-1)}^r| & \forall t \in T, k \in K, r \in S: t >1 & \label{eq:modelczr1} \\
&\underline{l}_{it}^r \leq 2l_{it}^r \leq \bar{l}_{it}^r  & \forall t \in T, i \in \Omega_{t}, r \in R & \label{eq:modelcsidelen} \\
&d_{ijt}^r \geq |g_{it}^r-g_{jt}^r|  & \forall t \in  T, \{i,j\} \in \Omega_{t}, r \in R: i \neq j & \label{eq:modeldis} \\
&c_{it}^r-l_{it}^r \leq g_{it}^r \leq c_{it}^r+l_{it}^r & \forall t \in T, i \in \Omega_{t},r \in R& \label{eq:modelcio} \\
&c_{it}^x-L^x(1-b_{ikt}+\beta_k) \leq g_{it}^x \leq c_{it}^x+L^x(1-b_{ikt}+\beta_k) & \forall t \in T,  k \in K, i \in \Omega_{t}& \label{eq:modelciox} \\
&c_{it}^y-L^y(2-b_{ikt}-\beta_k) \leq g_{it}^y \leq c_{it}^y+L^y(2-b_{ikt}-\beta_k) & \forall t \in T, k \in K, i \in \Omega_{t}& \label{eq:modelcioy} \\
&Mv_{it} \geq |c_{it}^r-c_{i(t-1)}^r| & \forall t \in T, i \in \Omega_{t},r \in R : t > 1& \label{eq:model} \\
&Mv_{it} \geq |l_{it}^r-l_{i(t-1)}^r| & \forall t \in T, i \in \Omega_{t},r \in R : t > 1 & \label{eq:model} \\
&u_{it}^r \geq |c_{it}^r-c_{i(t-1)}^r| & \forall t \in T, i \in \Omega_{t},r \in R : t > 1& \label{eq:model} \\
&a_{it}l_{it}^{x}+ 4\bar{x}_{ipt}^{2}l_{it}^{y} \geq 2a_{it}\bar{x}_{ipt}  & \qquad \forall t \in T, i \in \Omega_{t}, p=1,\dots,\Delta & \label{eq:modelcarea} 
\end{align}
\caption{The mixed-integer programming model of the Problem ZDFLP.}
\label{model}
\end{figure*}

The objective function of the ZDFLP aims to minimize the total cost of material handling and relayout during the planning horizon. The first part of the objective function (i.e., \eqref{eq:modelo1}) represents the total material handling cost. Equation \eqref{eq:modelo2} is the total relayout cost of departments with two components: (i) a variable relayout cost based on how much the center of a department is moved between two consecutive periods and (ii) a fixed cost independent of the distance the department is relocated. Equation \eqref{eq:modelo3} represents the cost of changing the boundaries of the zones between two consecutive periods. Moving zone boundaries from one period to the next may require modifications in the material handling system or aisle structures. In fact, such structural changes can be more costly than relocating departments and should be considered in the DFLP. 

Constraints \eqref{eq:modelczone}-\eqref{eq:modelcz14} are used to arrange the zone locations in the facility so they do not overlap. Constraints \eqref{eq:modelczone} ensure that zones $k$ and $h$ do not overlap at least in one axis direction. Constraints \eqref{eq:modelcz7}-\eqref{eq:modelcz10} set the zone boundaries based on the zone precedence relations defined in \eqref{eq:modelczone}, and \eqref{eq:modelcz11}-\eqref{eq:modelcz14} make sure that the zones are located with the facility. 

Constraints \eqref{eq:modelc1}-\eqref{eq:modelc4} are used to prevent the departments from overlapping in the axis direction in which they are arranged within a zone. A zone can be either $x$-axis or $y$-axis oriented. If two departments $i$ and $j$ are in the same $x$-axis oriented (or $y$-axis oriented) zone $k$, then these constraints are equal to $z_{ijt}^x+z_{jit}^x=1$ (or $z_{ijt}^y+z_{jit}^y=1$), forcing the departments to be non-overlapping due to constraints \eqref{eq:modelc4}. Note that these constraints are inactive for a department pair $i$ and $j$ assigned to different zones.   
Constraints \eqref{eq:modelcz1}-\eqref{eq:modelcz4} make sure that the departments are located within the boundaries of the zones to which they are assigned.   

Constraints \eqref{eq:modelczr1} capture whether the boundaries of the zones are changed between two consecutive periods. Relocating zone boundaries from one period to another may require significant changes in the material handling system or aisle structure. The cost of such structural changes can be more than relocating departments and should be considered in the DFLP. Constraints \eqref{eq:modelcsidelen} control the shape of the departments. Constraints \eqref{eq:modeldis} calculate the rectilinear distances between the I/O points of the departments.    

Constraints \eqref{eq:modelcio} make sure that the I/O point of a department is located within the boundaries of the department. Constraints \eqref{eq:modelciox} set the location of the I/O point to the center of a department in the $x$-axis direction (i.e., $g_{it}^x=c_{it}^x)$ if the department is located in an $x$-axis oriented zone. Similarly, constraints \eqref{eq:modelcioy} set $g_{it}^y=c_{it}^y$ if the department is located in a $y$-axis oriented zone. Note that these constraints do not directly require that departments have their I/O points adjacent to zone boundaries. Since the distances among the I/O points should be minimized to reduce the material handling cost for a given block layout, the I/O points tend to be located in the perimeters of departments. In addition, the departments are arranged in either $x$ or $y$-axis directions within the zones, ensuring that all departments are adjacent to a zone boundary. Thereby, the I/O points are located on zone boundaries to minimize the distances among the departments. Considering the I/O points of the departments in the ZDFLP provide more accurate modeling of material movements and interdepartmental distances within the facility, leading to a block layout that is quite different from the one designed based on center-to-center distances. Furthermore, the zone-based layout provides an appropriate structure for incorporating I/O points into a block layout design because zone boundaries can be used as candidate locations for aisles.    
Constraints \eqref{eq:modelcarea} model the department area requirements using the polyhedral outer approximation method of Sherali et al. \cite{sherali2003} based on $\Delta$ support points.  

\begin{figure*}
\centerline{\includegraphics[width=6in]{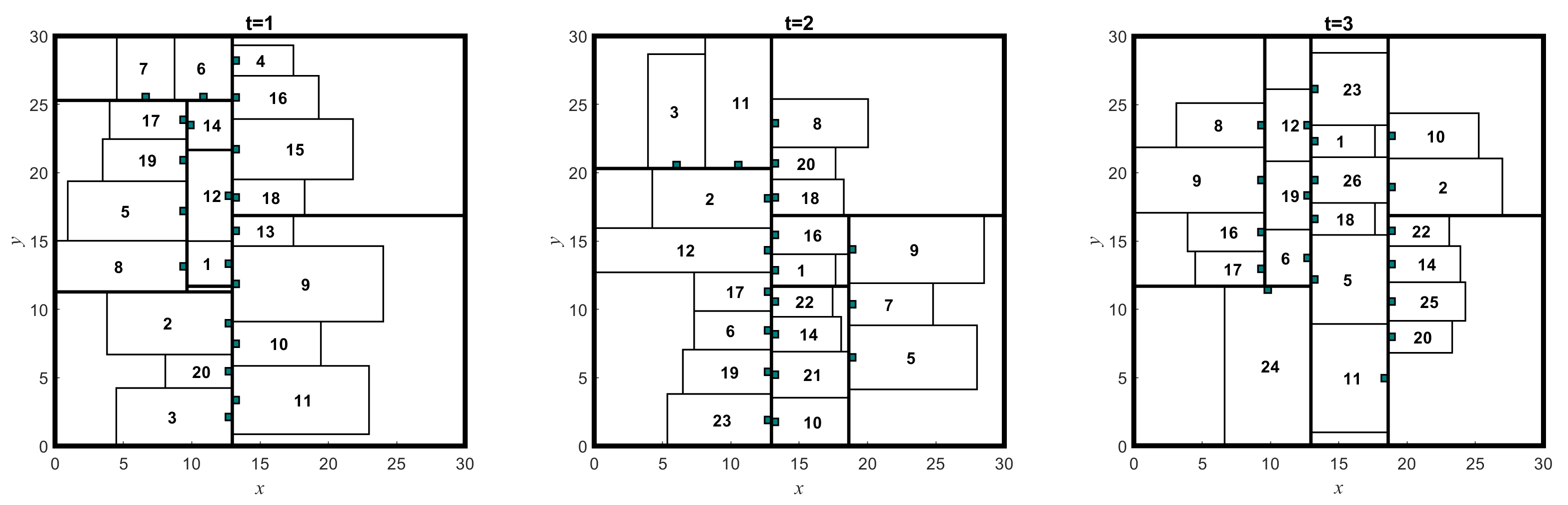}}
\caption{A solution found for DFLP 20-3 within 3600 CPU seconds ($t$=2: $21 \rightarrow 4$, $22 \rightarrow 15$, $23 \rightarrow 13$; $t$=3: $26 \rightarrow 21$, $25 \rightarrow 3$, $24 \rightarrow 7$, where, e.g., $21 \rightarrow 4$ means that department 21 replaces department 4 in period 2.)}
\label{lac20_3_1}
\end{figure*}

\section{Solution Methodology}
We developed a matheuristic based on variable neighborhood search (VNS) concepts to solve the ZDFLP with two phases. The matheuristic directly operates on Problem ZDFLP without requiring an encoding scheme of the problem and uses mixed-integer programming (MIP) to find new solutions. Therefore, this matheuristic will be referred to as MIP-VNS. In Phase I, an initial incumbent solution $S$ is found by solving Problem ZDFLP in CPLEX for a given time period. Incumbent solution $S$ is improved iteratively in Phase II by solving the Problem ZDFLP, but this time for only a subset $U$ of decision variables $z_{ijt}^r$, $b_{ikt}$, and $l_{it}^r$. Algorithm \label{NewCandidate} illustrates the procedure for generating a candidate solution $S_U$ from incumbent solution $S$. In Problem ZDFLP, all decision variables $z_{ijt}^r$, $b_{ikt}$, and $l_{it}^r$, excluding the ones in a given set $U$, are fixed to their corresponding values in $S$ (denoted by $z_{ijt}^r(S)$, $b_{ikt}(S)$, and $l_{it}^r(S)$). This reduced problem is denoted as Problem ZDFLP($U$) and can be solved quickly to find the optimal values of the decision variables in set $U$ using CPLEX, yielding a new candidate solution. To reduce the CPU requirement in Phase II, the CPU time to solve Problem ZDFLP($U$) can be limited by an upper bound.     

\begin{algorithm}\label{NewCandidate}
\SetAlgoLined
\KwIn{$S,U$}
\KwOut{$S^*,S_U$}
Fix $z_{ijt}^r \leftarrow z_{ijt}^r(S)$, $l_{it}^r \leftarrow l_{it}^r(S)$, $b_{ikt} \leftarrow b_{ikt}(S)$  \;  
Unfix all decision variables in set $U$ \;
Solve the ZDFLP optimally to find $S_U$ \;  
\caption{GenerateCandidate($S,U$)}
\end{algorithm}

A candidate solution is a local optimum within the incumbent solution $S$'s neighborhood defined by the variables in set $U$. In other words, the candidate solution is expected to improve upon or have the same objective value as the incumbent solution. Therefore, the Phase II search improves the incumbent solution quickly, but the progress stalls for a given neighborhood structure of $U$. When the search stalls, the neighborhood structure to select the decision variables to relax is changed. In Phase II, MIP-VNS systematically changes $U$ when the incumbent solution cannot be improved within the current neighborhood. The four neighborhood structures used in MIP-VNS are given in  \eqref{eq:u1} to \eqref{eq:u4}. For example, in neighborhood $U_1$, MIP-VNS updates incumbent solution $S$ by solving Problem ZDFLP($U$) for a randomly selected department $i^\prime$ and a randomly selected period $t^\prime$. In the largest neighborhood $U_4$, MIP-VNS tries to improve incumbent solution $S$ by solving Problem ZDFLP($U$) for two randomly selected departments $i^\prime$ and  $i^{\prime\prime}$ in two randomly selected periods $t^\prime$ and $t^{\prime\prime}$. If the total cost ($TC$) of the incumbent solution $S$ does not improve after trying all possible searches in a neighborhood structure of $U_k$, then MIP-VNS moves to neighborhood structure $U_{k+1}$. The MIP-VNS stops after trying $g_{max}$ neighborhood structures.  

\begin{multline}
U_1=\{z_{i^{\prime}jt^{\prime}}^r, z_{ji^{\prime}t^{\prime}}^r, b_{i^{\prime}kt^{\prime}}, l_{i^{\prime}t^{\prime}}^r : \forall j \neq i^{\prime}, k, r \}  \label{eq:u1}
\end{multline}
\begin{multline}
U_2=\{z_{i^{\prime}jt^{\prime}}^r, z_{ji^{\prime}t^{\prime}}^r, z_{i^{\prime}j(t^{\prime}+1)}^r, z_{ji^{\prime}(t^{\prime}+1)}^r, \\ b_{i^{\prime}kt^{\prime}},  b_{i^{\prime}k(t^{\prime}+1)}, l_{i^{\prime}t^{\prime}}^r, l_{i^{\prime}(t^{\prime}+1)}^r : \forall j \neq i^{\prime}, k, r \} \label{eq:u2}
\end{multline}

\begin{multline}
U_3=\{z_{i^{\prime}jt^{\prime}}^r, z_{ji^{\prime}t^{\prime}}^r, z_{i^{\prime\prime}jt^{\prime}}^r, z_{ji^{\prime\prime}t^{\prime}}^r, \\ b_{i^{\prime}kt^{\prime}}, b_{i^{\prime\prime}kt^{\prime}}, l_{i^{\prime}t^{\prime}}^r, l_{i^{\prime\prime}t^{\prime}}^r : \forall j \neq i^{\prime}, k, r \} \label{eq:u3}
\end{multline}

\begin{multline}
U_4=\{z_{i^{\prime}jt^{\prime}}^r, z_{ji^{\prime}t^{\prime}}^r,z_{i^{\prime}jt^{\prime\prime}}^r, z_{ji^{\prime}t^{\prime\prime}}^r,z_{i^{\prime\prime}jt^{\prime}}^r, z_{ji^{\prime\prime}t^{\prime}}^r,z_{i^{\prime\prime}jt^{\prime\prime}}^r, \\ z_{ji^{\prime\prime}t^{\prime\prime}}^r, b_{i^{\prime}kt^{\prime}}, b_{i^{\prime}kt^{\prime\prime}},b_{i^{\prime\prime}kt^{\prime}}, b_{i^{\prime\prime}kt^{\prime\prime}}, \\ l_{i^{\prime}t^{\prime}}^r, l_{i^{\prime}t^{\prime\prime}}^r, l_{i^{\prime\prime}t^{\prime}}^r, l_{i^{\prime\prime}t^{\prime\prime}}^r  : \forall j \neq i^{\prime}, j \neq i^{\prime\prime}, k, r \} \label{eq:u4}
\end{multline}

\begin{algorithm}\label{NewCandidate}
\SetAlgoLined
\KwIn{$S,U$}
\KwOut{$S^*,S_U$}
Solve the Problem ZDFLP until $\kappa$ feasible solutions are found. \;
$K \leftarrow 1 $ \;
$S^* \leftarrow S $\;
\For{$g = 1, \dots, g_{max}$}{ 
$noupdate \leftarrow noupdate+1$ \;
$A \leftarrow \{ (i, t):t \in T, i \in \Omega_t \} $ \;
\While{$A \neq \varnothing $ }{
Randomly and uniformly select $(i^\prime, t^\prime)$ from $A$ \;
$A  \leftarrow A -(i^\prime, t^\prime)$ \;
\If{$K=1$}{
$S_U$=GenerateCandidate($S,U_1$) \;
} 
\If{$K=2$}{
$S_U$=GenerateCandidate($S,U_2$) \;
} 

\If{$K=3$}{
Randomly select department $i^{\prime \prime}$ such that $i^{\prime \prime} \neq i^{\prime}$ \;
$S_U$=GenerateCandidate($S,U_3$) \;
} 
\If{$K=4$}{
Randomly select department $i^{\prime \prime}$ and period $t^{\prime \prime}$ such that $i^{\prime \prime} \neq i^{\prime}$ and  $t^{\prime \prime} \neq t^{\prime}$\;
$S_U$=GenerateCandidate($S,U_4$) \;
} 
\If{$TC(S_U) <TC(S^*)$ }{
		    $S \leftarrow S_U$ \;
		    $noupdate \leftarrow 0$ \;
} 
\If{$TC(S_U) <TC(S)$ }{
		    $S \leftarrow S_U$ \;
} 

}
\If{$noupdate>0$ }{
		    $K \leftarrow mod(K+1, 4)+1$ \;
		    $S \leftarrow S^*$ \;
} 

}
\caption{Finding}
\end{algorithm}

\section{Computational Experiments}
This paper introduces the ZDFLP with I/O points and flexible zone structures for the first time, so it's impossible to compare the results with previous studies directly. In particular, the literature lacks test instances for the DFLP with I/O points. Benchmarking was done using the test problems and results provided by \cite{Kulturel-Konak2019}, which used vertical and horizontal bands to form zones and I/O points in some test problems. Test problems DFLP 12-3c, DFLP 12-5c, and DFLP 20-3c were originally from \cite{Lacksonen1997}. Kulturel-Konak \cite{Kulturel-Konak2019} studied these problems considering band reallocation cost and I/O points that can be placed anywhere within the facility. Therefore, these problems provide some form of benchmark for the proposed approach in this paper. Problems FBS-DFLP-3b, and FBS-DFLP-4b were from \cite{mazinani2013}, which used the flexible bay structure with only vertical or horizontal bays. Kulturel-Konak \cite{Kulturel-Konak2019} also studied these problems using the flexible bay structure, including the cost of moving bay boundaries. In this paper,  $B_{kt}$  was set to 1/4 of the band allocation costs given in \cite{Kulturel-Konak2019} because the cost of changing each side of a zone is considered independently in the objective function. All other problem inputs and data were identical to \cite{Kulturel-Konak2019}. 

Table \ref{results} compares the best and average results found for the test problems given in \cite{Kulturel-Konak2019}. The results were found in five random replications of the VNS-MIP using $g_{max}$=50 on a Mac computer with an 8-Core Intel Core i9 CPU (2.4 GHz) and 32 GB system memory. As seen in Table \ref{results}, VNS-MIP found better solutions than the previously reported best solutions. In particular, the earlier best solutions of DFLP 12-3c, DFLP 12-5c, and DFLP 20-3c were also found by considering the distances among the I/O points of the departments. In fact, the results in \cite{Kulturel-Konak2019} were anticipated to be better than the ones in this paper since \cite{Kulturel-Konak2019} used I/O points that could be anywhere within departments. Therefore, the improvement in these problems could be attributed to the flexible zone structure used in this paper compared to band-based zones used in \cite{Kulturel-Konak2019}. Since FBS-DFLP-3b and FBS-DFLP-4b were not studied by considering I/O points before, the significant improvements on FBS-DFLP-3b and FBS-DFLP-4b were possibly due to the use of I/O points in addition to the flexible zone structure.   

\begin{table*}[htbp]
\caption{Results for Various Test Problems}
\begin{center}
\begin{tabular}{lrrrr}
\hline
 &\textbf{Previous}                                   & \textbf{Best{     }}       & \textbf{Average}  & \textbf{Average}      \\
\textbf{Problem}&\textbf{Best} \cite{Kulturel-Konak2019} & \textbf{Solution}   & \textbf{Solution}  & \textbf{CPU Sec.}      \\
\hline
FBS-DFLP-3b & 22,029.28                          & 20,929.59          & 21,609.27        & 908.23        \\
FBS-DFLP-4b & 41,885.26                          & 32,896.64 & 35,245.57 & 2,316.00 \\
DFLP 12-3c  & 4,735.06                           & 4,386.74  & 4,598.78  & 850.00 \\
DFLP 12-5c  & 7,978.32                           & 8,204.19   & 8,596.09  & 1,205.00        \\
DFLP 20-3c  & 9,700.21                           & 9,534.97   & 10,037.91        & 2,343.00     \\  
\hline
\end{tabular}
\label{results}
\end{center}
\end{table*}

\begin{figure}
\centerline{\includegraphics[width=2in]{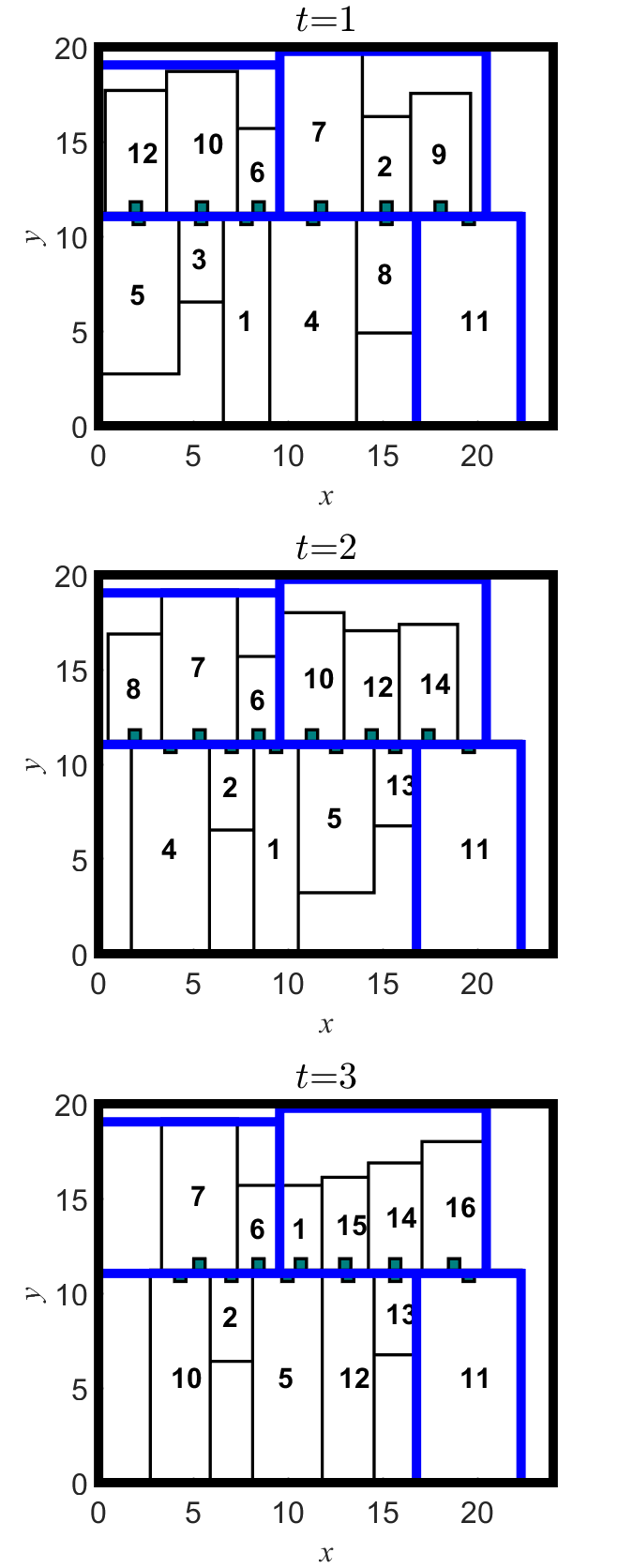}}
\caption{A solution found for DFLP 12-3b ($TC$=4,379.77) ($t$=2: $13 \rightarrow 3$, $14 \rightarrow 9$; $t$=3: $15 \rightarrow 4$, $16\rightarrow8$; where, e.g., $13 \rightarrow 3$ means that department 13 replaces department 3 in period 2.)}
\label{lac_97_12_3}
\end{figure}

\begin{figure}
\centerline{\includegraphics[width=3.2in]{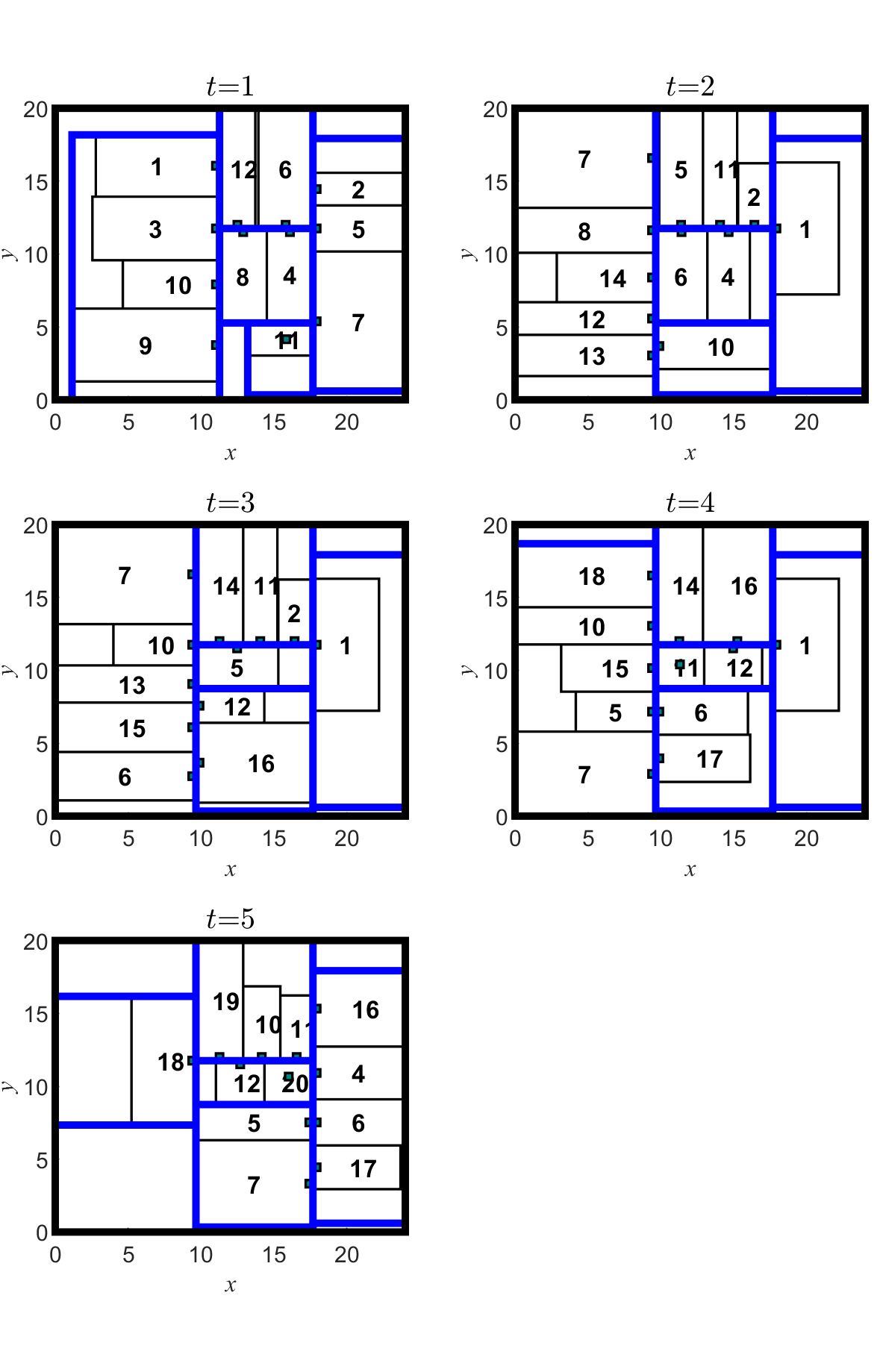}}
\caption{A solution found for DFLP 12-5b ($TC$=7,474.24) ($t$=2: $13 \rightarrow 3$, $14 \rightarrow 9$; $t$=3: $15 \rightarrow 4$, $16 \rightarrow 8$; 
$t$=4: $17 \rightarrow 2$, $18 \rightarrow 13$; $t$=5: $19 \rightarrow 1$, $20 \rightarrow 14$ where, e.g., $13 \rightarrow 3$ means that department 13 replaces department 3 in period 2.)}
\label{lac12_5_1}
\end{figure}

The primary objective of the model proposed in this paper is to represent DLFP as close to practice as possible. Fig. \ref{lac_97_12_3} illustrates the best solution found for DFLP 12-3b in five replications. In this problem, the four zones formed a spine layout where the I/O points of the departments were located at a center aisle between the upper and lower two zones. For this layout, the relative locations of the zones were partially defined (i.e., $\gamma_{12}^x=1$,$\gamma_{43}^x=1$,$\gamma_{14}^y=1$, and $\gamma_{23}^x=1$) prior to the optimization. However, the location of the zones and zone types were determined by the VNS-MIP. In this solution, the zone locations were fixed through the planning horizon despite the addition of new departments and changes in the sizes of the departments. Although four zones were used in this problem, the same layout could be achieved by two horizontal zones because the departments are not expected to fill zones completely. This example suggests that it would be beneficial to run the VNS-MIP with different numbers of zones or zone configurations to discover alternative facility designs. Fig. \ref{lac12_5_1} presents a solution found for DFLP 12-5b where the layout structure changes significantly from one period to the next one.  A unique aspect of this solution is that both vertical and horizontal aisle structures were utilized. It is clear in Fig. \ref{lac12_5_1} that the found block layout provides a good basis for possible aisles. Note that the I/O point of departments 11 and 20 were not on the boundaries of the zones. As discussed previously, this is a drawback of the approach used in determining I/O point locations. However, the proposed modeling approach is computationally efficient as it does not require any binary variables to model I/O point locations and provides a sufficient approximation to actual distances. As shown in Figs. \ref{lac20_3_1}, \ref{lac_97_12_3}, and \ref{lac12_5_1}, the proposed model can yield block layouts that are quickly transferable to detailed layouts. In addition, the ZDFLP model allows practitioners to define a rough layout plan that is appropriate to their processes and manufacturing systems.           

\section{Contributions and Conclusions}
This paper describes a MIP model for solving the dynamic facility layout problem. In the model, the departments are grouped into flexible zones, and the size and shape of these zones and departments are decision variables. The zones are separated by boundaries that could be used as locations for aisles and can be easily adjusted to accommodate changes in the material handling system. Departments are assumed to have their I/O points adjacent to zone boundaries. This zone-based approach allows designers to control the layout of the facility based on required manufacturing or service processes. The model also takes into account the cost of changing the overall layout structure over multiple planning periods. With these aspects and contributions, the proposed approach represents a more practical model of the DFLP. This paper offers the following insights:
\begin{enumerate}
\item The fact that facility layout significantly impacts manufacturing and service systems’ operation costs and efficiency, this paper will directly support material handling practitioners and researchers.
\item A main contribution of the paper is that the proposed mathematical model considers the cost of structural modifications in the facility and material handling system due to rearranging departments from one period to the next in addition to department relocation costs. This leads to dynamic block layouts that are easier to implement in practice than those obtained by unrestricted general mathematical models. Integration of departments' I/O points also leads block layouts closer to actual detailed layout implementation. 
\item The paper presents how a matheuristic can solve complex non-linear mixed integer programming FLPs by partitioning problems into smaller ones that can be solved optimally. Unlike heuristics approaches that require encoding the problem space into a problem representation schema, matheuristics directly operate on the decision variables of the problem. Therefore, matheuristics are highly suitable to solve DFLPs involving multi-level decisions such as the overall structure of the layout, shapes and locations of departments, and the I/O points. However, matheuristics are also subjective to premature convergence because of their strong local search characteristics. 
\end{enumerate}

Therefore, the model, proposed matheuristic and above strategies can be applied to other facility layout and material handling problems. Future research might apply the proposed model and approach to other industrial sectors, such as hospital and retail facility layouts. Another research direction is to refine the problem formulation to reduce the gap between the lower and upper bounds of solutions during optimization using techniques such as adding valid inequalities, improving the linear relaxation of the problem, or using stronger bounds on the constraints and variables. The proposed matheuristic may benefit from tightening the model. Finally, strategies to improve the performance of the proposed matheuristic is an interesting further research topic.

\bibliographystyle{IEEEtran}
\bibliography{imhrc}

\end{document}